\documentclass[twocolumn]{article}
\usepackage{cite}
\usepackage{amsmath,amssymb,amsfonts}
\usepackage{amsthm}

\usepackage{algorithm,algorithmic}
\usepackage{graphicx}
\usepackage{textcomp}
\usepackage{xcolor}
\usepackage{subfig}
\def\BibTeX{{\rm B\kern-.05em{\sc i\kern-.025em b}\kern-.08em
    T\kern-.1667em\lower.7ex\hbox{E}\kern-.125emX}}

\usepackage[english]{babel}
\usepackage[latin1]{inputenc}
    
\newcommand{\R}{\mathbb{R}}
\newcommand{\fun}{\mathcal{F}}
\newcommand{\xtrue}{\widetilde{x}}
\newcommand{\argmin}[1]{\underset{#1}{\mathrm{argmin\,}}}

\newtheorem{theorem}{Theorem}

\newtheorem{assumption}{Assumption}
\graphicspath{{./figures/}}

\begin{document}

\title{Non-convex approach to binary compressed sensing}

\author{Sophie M. Fosson}
\maketitle

\begin{abstract}
We propose a new approach for the recovery of binary signals in compressed sensing, based on the local minimization of a non-convex cost functional. The desired signal is proved to be a local minimum of the functional under mild conditions on the sensing matrix and on the number of measurements. We develop a procedure to achieve the desired local minimum, and, finally, we propose  numerical experiments that show the improvement obtained by the proposed approach with respect to classical convex methods.
\end{abstract}
%
%
\section{Introduction}\label{sec:intro}
Binary compressed sensing (BCS, \cite{nak12}) refers to compressed sensing (CS, \cite{don06}) in the case of sparse binary signals. We can mathematically formulate BCS as follows: recover the sparsest $x\in\{0,1\}^n$ from $y=Ax$, given $A\in\R^{m,n}$ with $m<n$. This is the simplest formulation of the most general problem of CS for finite-valued or discrete signals, that is, signals whose components are known to belong to a finite or discrete alphabet. The first theoretical analyses on this topic are very recent \cite{kei17}, despite the relevance in practical applications is widespread, e.g., in digital communications, wireless communications, sensor networks, digital image processing, spectrum sensing, localization, and quantized systems. 

It is worth noticing that the classical recovery algorithms for CS can be applied to finite-valued signals as well, but they are not always prone to embed prior information on the discrete nature of the signal, which is expected to improve the accuracy of the solution. On the other hand, the discrete nature might lead to combinatorial approaches, which turn out to be computationally burdensome. In this paper, we instead show that solving a non-convex problem over the convex hull of the alphabet can lead to the exact recovery. In \cite{kei17}, the problem is tackled by solving the basis pursuit \cite{fou13} on the convex hull of the alphabet, and theoretical results on phase transition and stability are provided under null space properties.

In this paper, we propose a new formulation of the problem, introducing a cost functional tailored for binary signals in $\{0,1\}$. For simplicity, we restrict the problem to the noise-free setting. Generalizations to large alphabets and  noisy settings are under study. We however notice that the binary case is itself relevant, and BCS has drawn some attention in the last years, see, e.g., \cite{sto10,nak12,tia09,shi15,ahn16, lee16,liu18}. In these papers, methods that take into account the prior information on the binary nature of the signal are proposed, and are shown to perform better than classical CS strategies in numerical experiments. Moreover, the problem of recovering binary signals is widespread in different frameworks, e.g., localization \cite{bay15}, hybrid systems \cite{fox13}, and jump linear systems \cite{fox15}, with possible applications to fault detection.

Our main contribution is the introduction and analysis of a cost functional suitable for BCS. In particular, we prove that a minimum  is exactly the desired signal. A drawback of our approach is the non-convexity of the functional; however, we develop a suitable iterative procedure, which can be used to look for the the desired minimum. This method is shown to improve the recovery with respect to classical CS convex approach.

The paper is organized as follows. In Section \ref{sec:ps}, we introduce the problem. In Section \ref{sec:theo} we prove our main theoretical guarantees. In Section \ref{sec:sims}, we describe the  recovery algorithm and show some simulation results. Finally, we draw some conclusions.

\section{Problem statement}\label{sec:ps}
In this work, we propose a new efficient approach to BCS. Specifically, we consider the following cost functional:
\begin{equation}\label{mcplasso}
\begin{split}
\fun&:[0,1]^n \to \R^+\\
\fun(x)&:=\frac{1}{2}\left\|y-A x \right\|_2^2+\lambda\sum_{i=1}^n\left(x_i-\frac{1}{2}x_i^2\right),~~~\lambda>0.
\end{split}
\end{equation}
$\fun(x)$ is similar to the popular Lasso functional \cite{tib96, fou13}, but it presents a concave penalty $g(x)= \sum_{i=1}^n\left(x_i-\frac{1}{2}x_i^2\right)$, which belongs to the family of minimax concave penalties (MCP).  In the last years, the MCP family  has been used and analyzed for sparse recovery and variable selection \cite{zha10MCP,zha12,woo16}, and for in-network recovery of jointly sparse signals  \cite{fox16}. 

Even though other concave penalties are more popular in the literature (for example, the $\log$ penalty \cite{can08rew}), $g(x)$ is preferable for its quadratic polynomial structure, which makes it mathematically more tractable. For example, it is possible to compute its global minimum by quadratic programming or other techniques of polynomial optimization \cite{lasbook}, even though this is computationally burdensome. This point is not examined in this paper, but may be considered for future developments. As we explain in the next section, in this work we rather focus on fast algorithms to find the desired local minima of $\fun$.

\section{Theoretical guarantees}\label{sec:theo}
In this section, we show that finding a (local) minimum of \eqref{mcplasso} can lead to the exact recovery of the unknown binary, sparse signal, under mild assumptions. Specifically, our assumptions are weaker than usual conditions required in CS: neither coherence conditions nor restricted isometry property nor null space property \cite{fou13} are involved in our analysis. 
The main results can be summarized as follows: the true $\xtrue\in\{0,1\}^n$ is a local minimum of \eqref{mcplasso}, and it is the unique signal in $\{0,1\}^n$ which is a local minimum of $\fun$ over the convex hull $[0,1]^n$. Therefore, using any descent algorithm that achieves a local minimum, if the achieved value is binary, then it is the exact signal.

The non-convex $\fun$ can have multiple minima, thus finding the right solution may be not easy. However, the following points should be considered.
\begin{itemize}
\item Fast algorithms are available to get local minima, that can be run multiple times starting from different initial points. This increases the chance of getting the right solution.
\item The non-convexity of $\fun$ is described by its Hessian matrix $A^T A -\lambda I$, which has $n-m$ negative eigenvalues (assuming $\lambda$ sufficiently small). Increasing the number of measurements $m$ we reduce the number of negative eigenvalues, thus we somehow mitigate the non-convexity of $\fun$ (which sometimes is even convex in $[0,1]^n$).
\end{itemize}
These observations suggest that suitable design of the system and implementation of the recovery method reinforce the possibility of exact reconstruction. In this paper, these points are verified via numerical simulations, while a theoretical investigation is left for future extended work.

Let us now prove the theoretical guarantees.

Let $\xtrue\in\{0,1\}^n$, with support $S$ and sparsity level $k$, and let $y=A\xtrue$. We indicate by $A_{[m]}$ any submatrix of $A$ obtained by selecting $m$ columns, while we write $Q\succ 0$ to say that matrix $Q$ is positive definite.  We recall that a family of vectors $\{v_1,\dots,v_n\}\subset \R^m$, $m<n$, are said to be in general position if any $\sum_i \sigma_i v_i$, with $\sigma_i\in\{0,\pm 1\}$, is different from $\pm v_j$, for any $j=1,\dots,n$ \cite{tib13}.

\begin{assumption}\label{ass1}$~$
\begin{itemize}
 \item[({\tt a})] $A_{[m]}^T A_{[m]}-\lambda I\succ 0$;
 \item[({\tt b})] the sparsity level is not larger than the number of measurements: $k\leq m$;
 \item[({\tt c})]the columns of $A$ are in general position.
\end{itemize}
\end{assumption}
Conditions ({\tt a}) and ({\tt c}) are usual: for example, any random matrix with entries generated according to a continuous distribution satisfy them with probability one \cite{tib13}. Instead, condition ({\tt b}) is a necessary  requirement to recover a sparse signal (in CS theory, at least $m\geq c k \log \left(\frac{n}{m}\right)$ for some $c>0$ is required \cite{fou13}).
\vskip0.5cm

\begin{theorem}\label{local}
Let $\xtrue\in\{0,1\}^n$, with support $S$ and sparsity level $k$, and let $y=A\xtrue$. Under assumptions \ref{ass1}.({\tt a})-({\tt b}), $\xtrue$ is an isolated local minimum of $\fun(x)$ as defined in \eqref{mcplasso}.
\end{theorem}

\begin{proof}
Let $\xtrue\in\{0,1\}^n$ and let $h\in\R^n$ be any small increment such that $\xtrue+h\in [0,1]^n$. Our goal is to show that $\fun(\xtrue+h)> \fun(\xtrue)$ for any $h\neq 0$. Let $S$ be the support of $\xtrue$, and $S^C=\{1,\dots,n\}\setminus S$.
 \begin{equation*}
  \begin{split}
   &\fun(\xtrue+h)-\fun(\xtrue) =\\
    & = \frac{1}{2}\left\|A h \right\|_2^2+\lambda \sum_{i=1}^n h_i-\frac{\lambda}{2}\left\|h\right\|_2^2-\lambda\langle \xtrue,h \rangle\\
    & = \frac{1}{2}\left\|A h \right\|_2^2+\lambda \sum_{i=1}^n h_i-\frac{\lambda}{2}\left\|h\right\|_2^2-\lambda\sum_{i\in S} h_i\\
    & = \frac{1}{2}\left\|A h \right\|_2^2+\lambda \sum_{i\in S^C} h_i-\frac{\lambda}{2}\left\|h\right\|_2^2.\\
  \end{split}
 \end{equation*}
In particular, let us notice that $h_i\geq 0$ when $i\in S^C$ (while $h_i\leq 0$ when $i\in S$).

We now prove that $f(h)=\frac{1}{2}\left\|A h \right\|_2^2+\lambda \sum_{i\in S^C} h_i-\frac{\lambda}{2}\left\|h\right\|_2^2> 0$ for any $h\neq 0$ with sufficiently small magnitude. This is sufficient to obtain the thesis. Let us say that $|h_i|<\epsilon$, where $\epsilon$ will be assessed in a while.

Since $f$ is differentiable, we compute its gradient to look for minima in the interior part of its domain. We easily obtain a unique stationary point  $-(A^T A-\lambda I)^{-1}\lambda v$ where $v\in\{0,1\}^n$ has entries equal to 1 in $S^C$ and zero elsewhere. Assuming that $\lambda$ is smaller than the positive eigenvalues of $A^T A$, the Hessian matrix $A^T A-\lambda I$ has both positive and negative eigenvalues (in particular, the negative ones are equal to $-\lambda$). Therefore, the stationary point is a saddle. In order to find minima, we then move on the boundaries. Let $A_{[m]}$ be any selection of $m$ columns of $A$. We assume that $A_{[m]}^T A_{[m]}-\lambda I \succ 0$ (this is always false for  any $A_{[j]}$, $j>m$) according to Assumption \ref{ass1}.({\tt a}). Therefore,  $n-m$ entries must be on boundaries ($\{0,\epsilon\}$ for $h_{S^C}$, $\{-\epsilon,0\}$ for $h_S$) to have a candidate minimum. Fixed, these $n-m$ entries, let $\Omega$ be the set of the remaining $m$ entries: we compute the gradient on $\Omega$ to look for the minimum. We obtain $h_{\Omega}=(A_{\Omega}^T A_{\Omega}-\lambda I)^{-1}(A_{\Omega}^T A_{\Omega^C}h_{\Omega^C}-\lambda v)$ where $v\in\{0,1\}^m$ is 1 on $S^C\cap \Omega$. By Assumption \ref{ass1}.({\tt a}), $A_{\Omega}^T A_{\Omega}-\lambda I \succ 0$, therefore this corresponds to a minimum. However, since the entries of $h_{\Omega^C}$ are in $\{0,\pm \epsilon\}$, we can choose an $\epsilon$ small enough so that the entries of $h_{\Omega}$ all have magnitude greater than $\epsilon$. More precisely, if $\epsilon$ is much smaller than $\lambda$, $h_{\Omega}$ can be approximated by $(A_{\Omega}^T A_{\Omega}-\lambda I)^{-1}(-\lambda v)=-(\lambda A_{\Omega}^T A_{\Omega}- I)^{-1}v$, which in turn can be approximated by $v$ is $\lambda$ is small. In conclusion the so-computed candidate minimum is outside the domain, and also the entries over $\Omega$ should be on the boundaries.
In conclusion, we are observing that a candidate minimum for $f$ is in $\{0,\pm\epsilon\}^n$.

At this point, we notice that if $h_{S^C}=0$, then $f(h)=\frac{1}{2}\left\|A_S h_S \right\|_2^2-\frac{\lambda}{2}\left\|h_S\right\|_2^2$, which has minimum for $h_S=0$. On the other hand, if there exists a $j\in S^C$ such that $h_j=\epsilon$, this is sufficient to have $f(h)> \lambda\left(\epsilon-\frac{k}{2}\epsilon^2\right)$, which is positive for any $\epsilon<\frac{2}{k}$ (which again can be assumed for the arbitrariness of $\epsilon$). This proves that $f(h)\geq 0$ for any sufficiently small $h$, and $f(h)=0$ if and only if $h=0$.

As a consequence, $\fun(\xtrue+h)> \fun(\xtrue)$ for any $h\neq 0$.
\end{proof}
\vskip0.5cm

\begin{theorem}\label{global}

Let $\xtrue\in\{0,1\}^n$, with support $S$ and sparsity level $k$, and let $y=A\xtrue$.
Under assumptions \ref{ass1}.({\tt a})-({\tt b})-({\tt c}), if $\lambda$ is sufficiently small, then $\xtrue$ is the global minimum of $\fun(x)$ defined in \eqref{mcplasso} over $\{0,1\}^n$.
\end{theorem}

\begin{proof}
Let us consider $f(h)$ as defined in the proof of Theorem \ref{local}. In order to investigate the minima of $\fun$ over $\{0,1\}^n$, we consider $h_{S}\in\{0,-1\}^k$ and $h_{S^C}\in\{0,1\}^{n-k}$.
Therefore, $f(h)\geq \frac{1}{2}\|Ah\|_2^2-\lambda\frac{k}{2}$. Now, for any $h\neq0$, if the columns of $A$ are in general position, then  $\frac{1}{2}\|Ah\|_2^2$ has a positive value. Therefore, we can always assume $\lambda < \frac{\|Ah\|_2^2}{k}$, which proves that $f(h)>0$. This yields to $\fun(z)>\fun(\xtrue)$ for any $z\in\{0,1\}^n$.
\end{proof}

The following theorem reinforces this result by proving that $\xtrue$ is the unique minimum over $\{0,1\}^n$. 
\vskip0.5cm
\begin{theorem}\label{check}
Let us consider any $z\in\{0,1\}^n$, $z\neq \xtrue$. Under assumptions \ref{ass1}.({\tt a})-({\tt b})-({\tt c}),  if $\lambda$ is sufficiently small, then, $z$ is not a minimum of $\fun$.
\end{theorem}
\begin{proof}
Let us consider any $z\in\{0,1\}^n$, $z\neq \xtrue$. We prove that there exists a direction $h$ along which $\fun(z+h)<\fun(z)$. First, we have:
\begin{equation}\label{eq:unique}
  \begin{split}
   &\fun(z+h)-\fun(z) =\\
    & = \frac{1}{2}\left\|A h \right\|_2^2+\langle Ah, A(z-\xtrue) \rangle +\lambda \sum_{i\in S^C} h_i-\frac{\lambda}{2}\left\|h\right\|_2^2.\\
  \end{split}
 \end{equation}
 Now, let us define $h=-\epsilon (z-\xtrue)$, where $\epsilon>0$. This $h$ is an admissible increment, since  if $z_i=1$, then $h\leq 0$, while if $z_i=0$, then $h\geq 0$; therefore, for $\epsilon\leq 1$, $z+h\in [0,1]^n$. Substituting this $h$ in \eqref{eq:unique}, we obtain:
 \begin{equation}\label{eq:conclude}
  \begin{split}
  &\fun(z+h)-\fun(z) \leq \frac{\epsilon^2}{2}\left\|A(z-\xtrue)  \right\|_2^2 -\epsilon \left\|A(z-\xtrue)  \right\|_2^2+\\& +\lambda \epsilon\|z-\xtrue\|_2^2-\frac{\lambda}{2}\epsilon^2 \left\|z-\xtrue\right\|_2^2=\\
  & = \epsilon \left(\frac{\epsilon^2}{2}-1\right) \left(\left\|A(z-\xtrue)  \right\|_2^2 -\lambda\|z-\xtrue\|_2^2 \right) 
  \end{split}
 \end{equation}
 where we use the fact that $\|\cdot\|_1$ is equal to $\|\cdot\|_2^2$ if the argument is binary.
 As the columns of $A$ are in general position, then $\left\|A(z-\xtrue)  \right\|_2^2$ is a positive value, and we can always choose $\lambda$ so that $\left\|A(z-\xtrue)  \right\|_2^2 -\lambda\|z-\xtrue\|_2^2>0$. Thus, the last expression in \eqref{eq:conclude} is negative for any $\epsilon<\sqrt{2}$. Since this is true for any arbitrarily small $\epsilon>0$, we conclude that $\fun(z+h)-\fun(z)<0$ if the direction of $h$ is $z-\xtrue$. This proves that $z$ cannot be a minimum.
 \end{proof}

\section{Numerical results}\label{sec:sims}
\begin{figure*}[ht]
\centering
\includegraphics[width=0.4\textwidth]{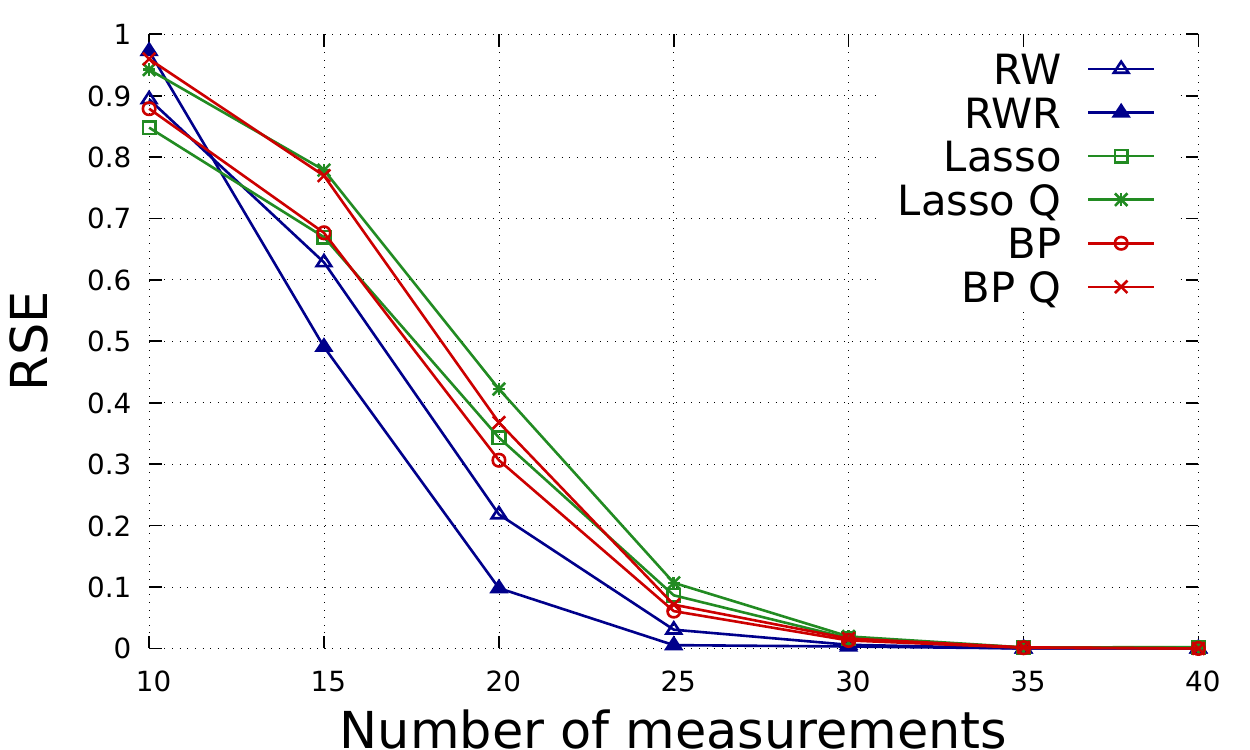}$~~$
\includegraphics[width=0.4\textwidth]{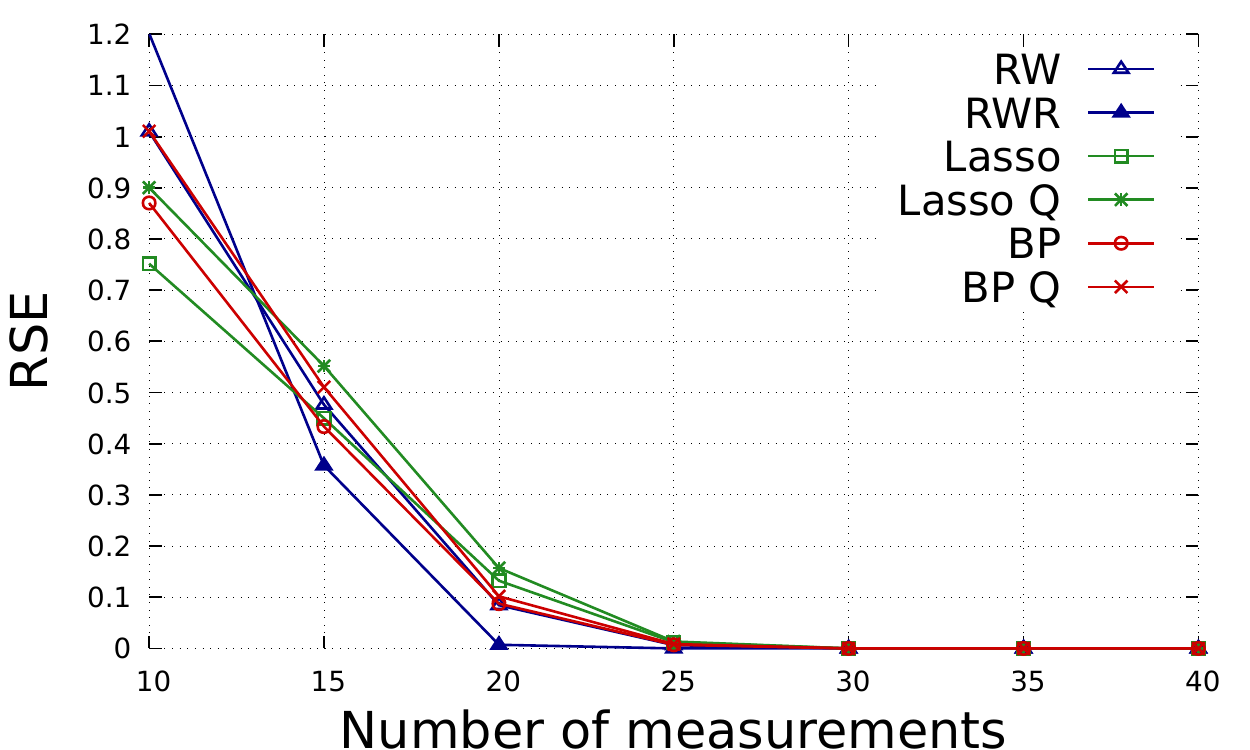}\\
\includegraphics[width=0.4\textwidth]{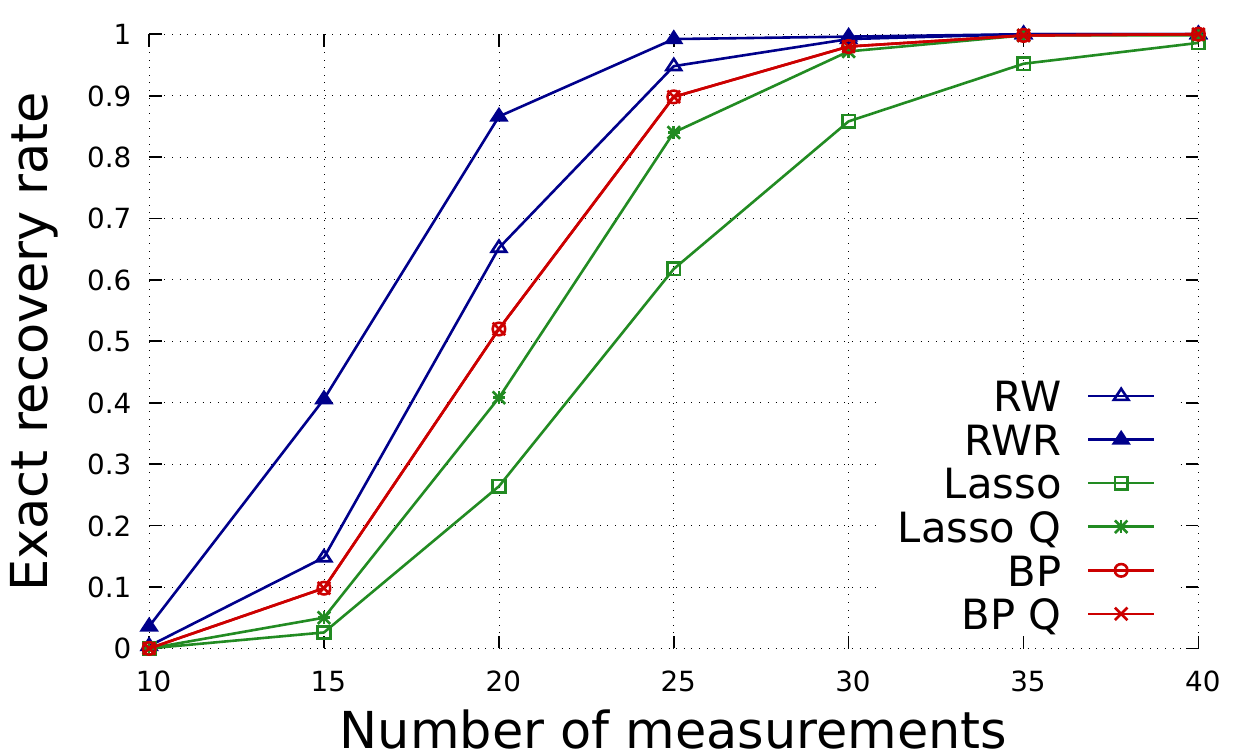}$~~$
\includegraphics[width=0.4\textwidth]{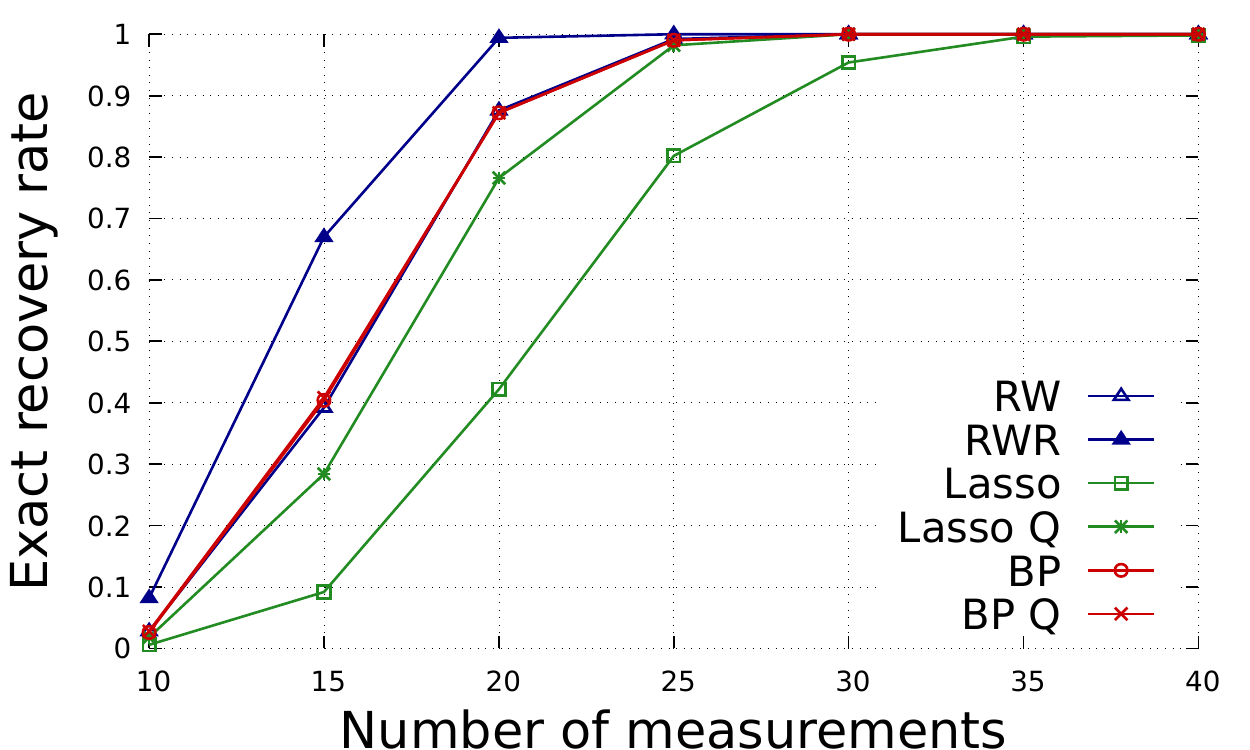}\\
\includegraphics[width=0.4\textwidth]{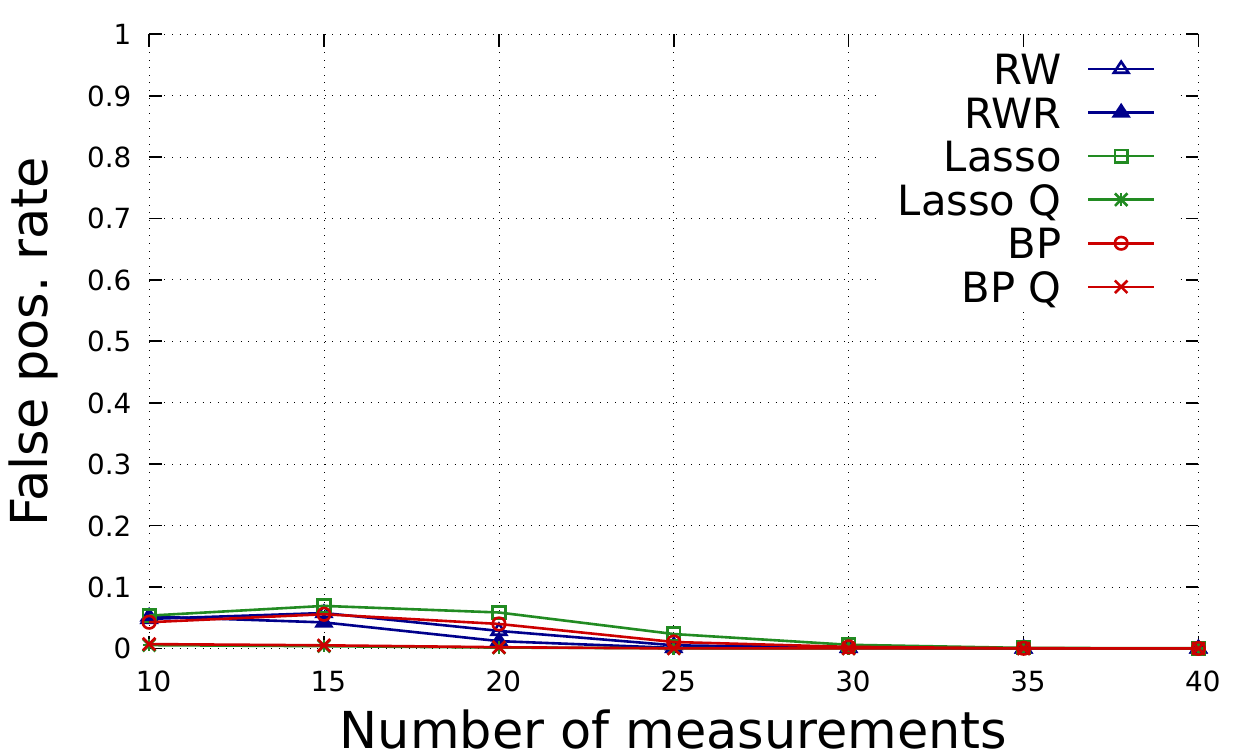}$~~$
\includegraphics[width=0.4\textwidth]{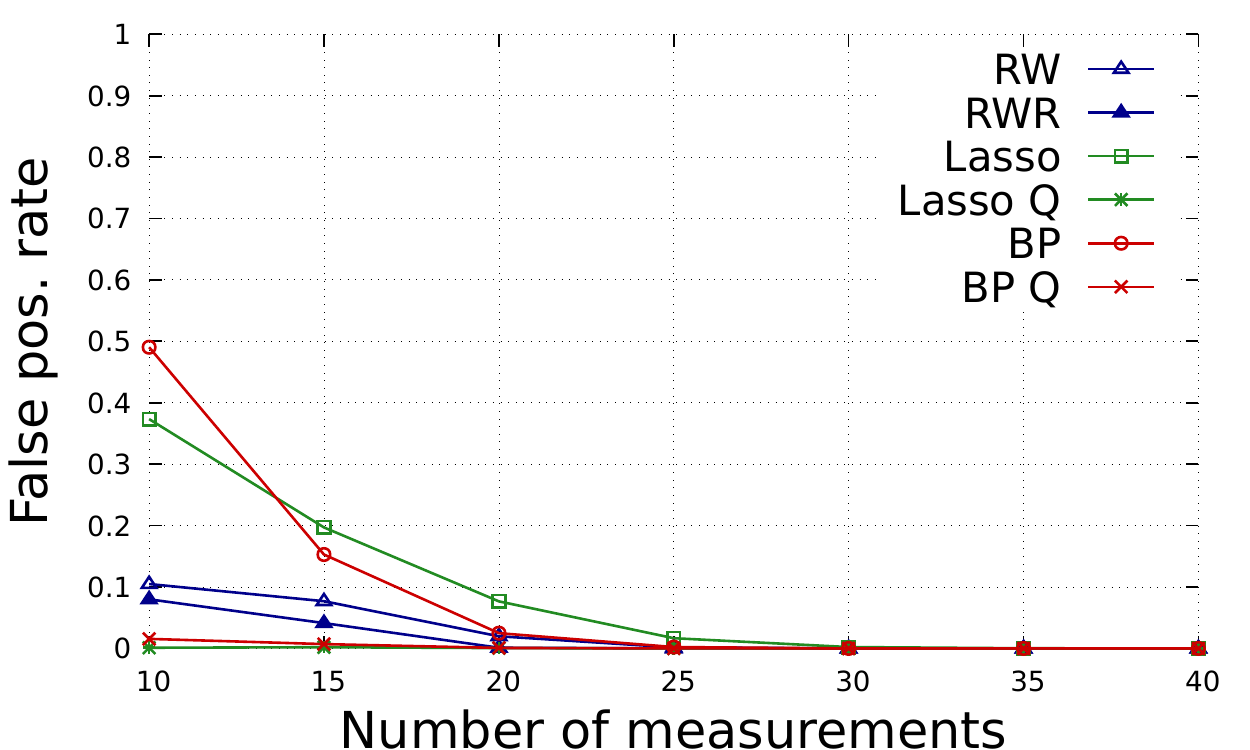}\\
\includegraphics[width=0.4\textwidth]{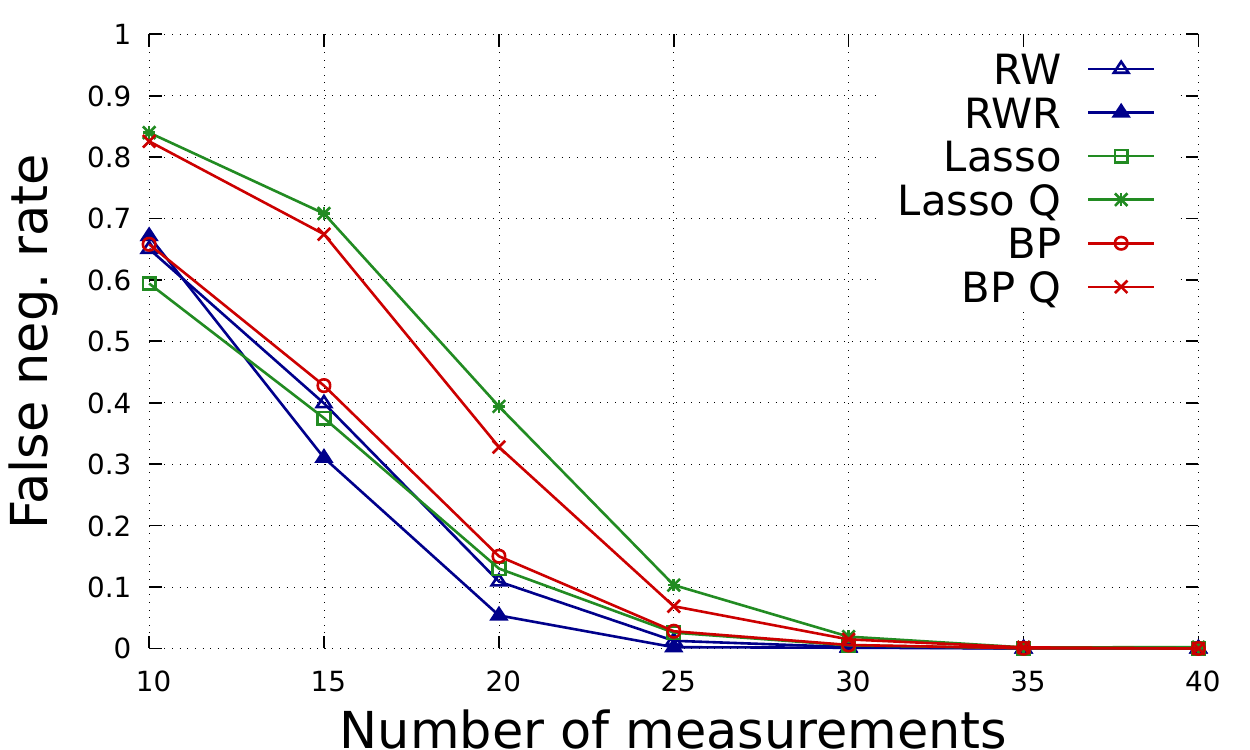}$~~$
\includegraphics[width=0.4\textwidth]{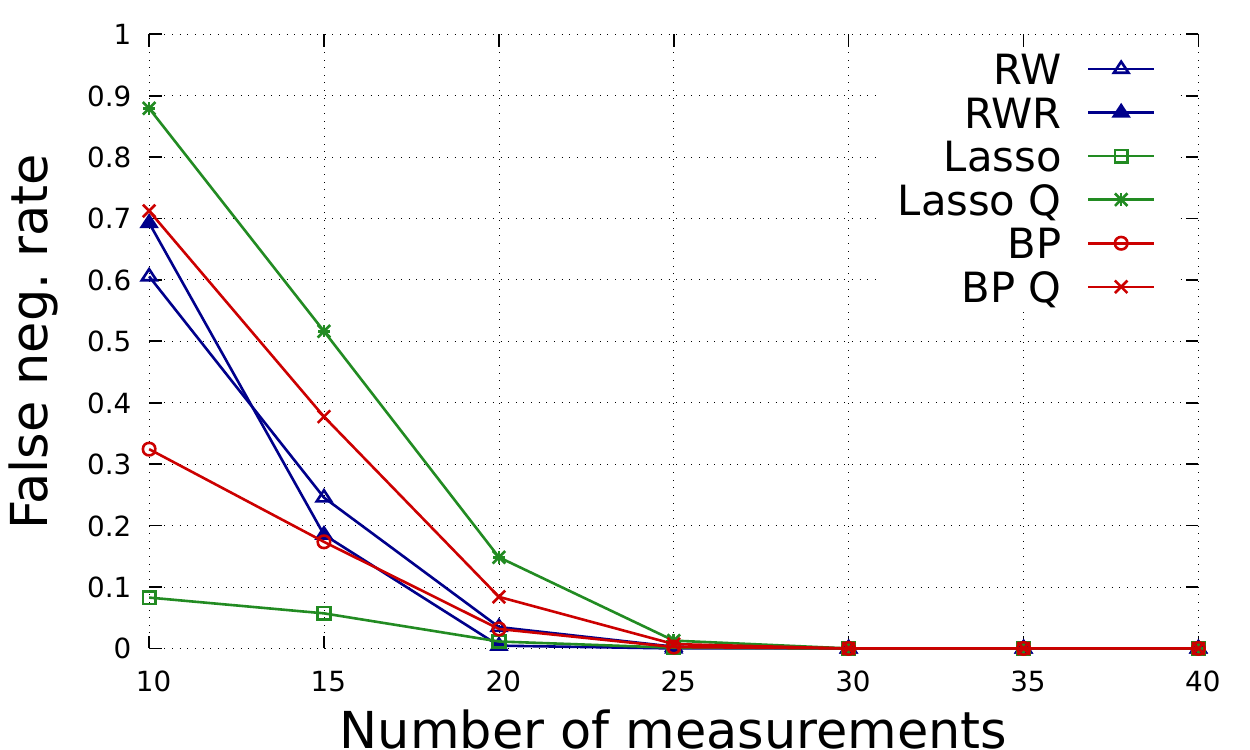}\\
\includegraphics[width=0.4\textwidth]{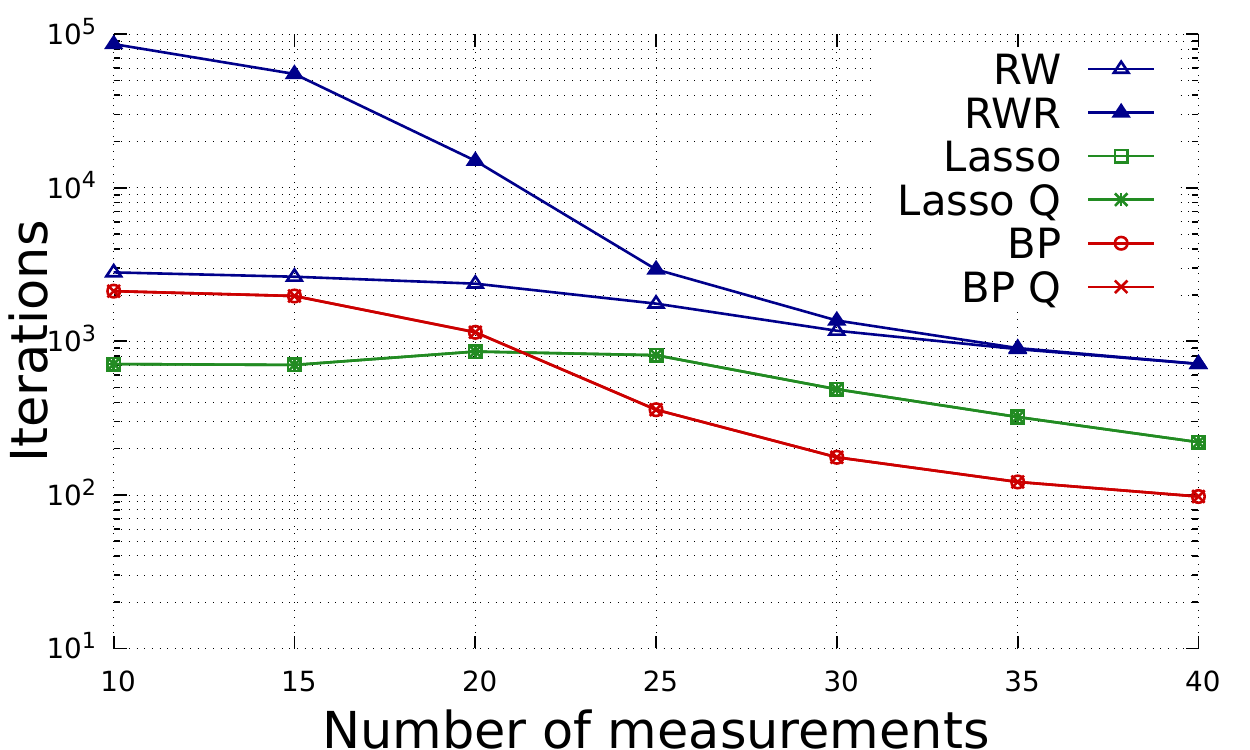}$~~$
\includegraphics[width=0.4\textwidth]{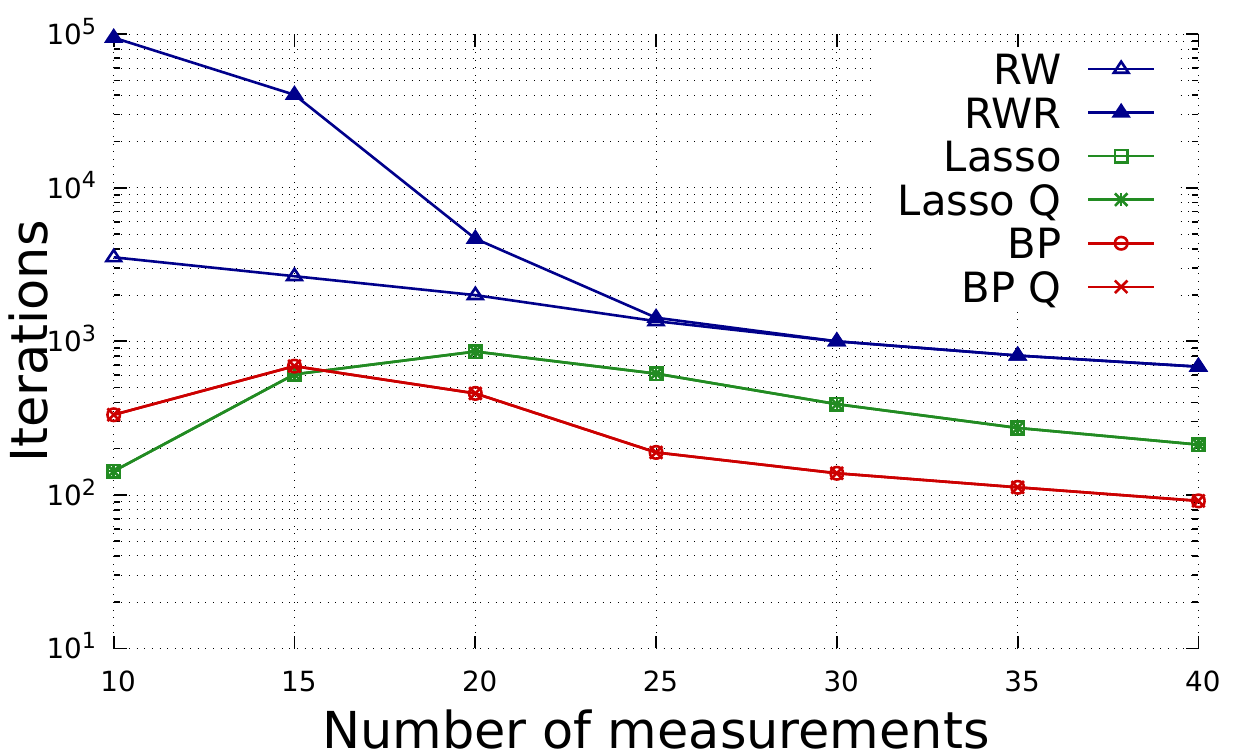}
\caption{Noise-free measurements, $n=100$, $k=5$; left column: unknown $k$, right column: known $k$. For Lasso, RW, RWR: $\lambda=10^{-2}$. Lasso Q and BP Q refer to the solutions obtained by quantizing the Lasso and BP solutions over $\{0,1\}^n$ (this is particularly useful for Lasso, which produces biased  solutions in the absence of noise).}\label{fig1}
\end{figure*}
In this section, we propose some numerical simulations that support the efficiency of our approach. 

\subsection{Algorithm}
An algorithm that can be used to compute a (local) minimum of $\fun$ in $[0,1]^n$ is the iterative reweighting  procedure (abbreviated as RW) introduced in Algorithm \ref{rw}.

 \begin{algorithm}
  \caption{Reweighting (RW)}\label{rw}
  \begin{algorithmic}[1] 
  \STATE  Initialize: $x(0)\in [0,1]^n$; $\lambda>0$
  \FORALL{$t=1,\dots,T_{stop}$}
  \STATE Update weights: $$w_i(t)=1- x_i(t) \text{ for any }i\in \{1,\dots, n\}$$
  \STATE Lasso: $$x(t+1)= \argmin{x\in [0,1]^n} \frac{1}{2}\left\|y-Ax\right\|_2^2 + \lambda\sum_{i=1}^n w_i(t) x_i$$
  \ENDFOR     
  \end{algorithmic}
\end{algorithm} 
$1- x_i(t+1)$ is the derivative of $g(x_i(t+1))=x_i(t+1)-\frac{1}{2} x_i(t+1)^2$. In the literature, such reweighting method has been studied for different concave penalties  (e.g., $\log(|x|+\epsilon)$, $(|x|+\epsilon)^q$, $q\in (0,1)$) and proved to reach a local minimum of concave-penalized functionals of kind $\frac{1}{2}\|Ax-y\|_2^2+\lambda \sum_i g(|x_i|)$ \cite{zou08_LLA, can08rew, faz03, fox18}. Even though limited to local minimization, practical experiments show that reweighting algorithms are generally very efficient and produces accurate estimates in many applications, e.g., magnetic resonance imaging, sensor selection and distributed CS in sensor networks \cite{can08rew, asi13, ahm15, cal16, fox14, fox16}.
 
For our experiments, we use reweighting, which turns out to be efficient. We however notice that other approaches might be tested, ranging from alternating minimization (subproblems of dimensions $m$ are convex) and block-coordinate descent to the alternating direction method of multipliers (ADMM, \cite{boy10}). We leave a comparison analysis for future work.
\subsection{Results}
We consider the following setting for our simulations\footnote{Code available at https://github.com/sophie27/Non-convex-approach-to-binary-compressed-sensing}. The desired signal is $\xtrue\in\{0,1\}^n$ with $n=100$ and sparsity level $k=5$. $m$ measurements are taken, $y=A\xtrue$, through a Gaussian sensing matrix $A\sim \mathcal{N} (0,\frac{1}{m})$, $m\in[10,40]$. The Lasso step in Algorithm \ref{rw} is iteratively solved via ADMM, which is stopped when the sum of squared  primal and dual residuals is below $10^{-6}$ \cite{boy10}. We study two possibilities: the first one, indicated as RW, is the $\ell_1$ reweighting  procedure introduced in Algorithm \ref{rw} with initial condition $x(0)=0$; in the second one, indicated by RWR, we re-run RW when the obtained solution is verified to be not the desired one via Theorem \ref{check}, restarting from $x(0)\in[0,1]^n$ generated uniformly at random. This is an heuristic method to test different initializations over $[0,1]^n$. We set $T_{stop}=4$ (as noticed in \cite{can08rew}, the most of benefit is generally obtained in the first reweigthing iterations), while the maximum number of re-initializations in RWR is fixed to 20.
We compare the proposed method to Lasso and Basis Pursuit (BP) \cite{fou13}, both solved via ADMM. BP is conceived for the free-noise case, while Lasso is known to have a bias. Lasso and BP do not envisage the prior information on the binary nature of the signal, but one can quantize the obtained  solution over $\{0,1\}$: we compute the performance metrics in both standard and quantized solutions.
Finally, we notice that, for all the methods (RW, RWR, Lasso, BP), when $k$ is known, we can add the equation $\sum_{i}x_i=k$ to the system $Ax=y$, which is expected to improve the performance.
The following performance metrics are evaluated (we indicate by $\widehat{x}$ the estimate of $\xtrue$):
\begin{itemize}
\item relative square error $RSE = \|\widehat{x}-\xtrue\|_2^2/\|\xtrue\|_2^2$;
\item false positive rate, that is, the normalized number of occurrences $\widehat{x}_i \neq 0$ when $\xtrue_i=0$;
 \item false negative rate, that is, the normalized number of occurrences $\widehat{x}_i = 0$ when $\xtrue_i \neq 0$;
 \item exact recovery, which is defined as $RSE<10^{-3}$ and no false positives/negatives;
 \item total number of ADMM iterations.
 \end{itemize}
The results, displayed in Figure \ref{fig1}, are averaged over 500 runs. When $k$ is unknown (left column), RW and RWR always achieve better performance than classical methods when $m\geq 15$. In particular, the random re-initialization in RWR gives a substantial gain, and always gets the exact solution at  $m=25$, where classical methods do not overpass $90\%$ of success. When $k$ is known (thus equation $\sum_i x_i=k$ is added), as expected, the general recovery accuracy is improved. In this case, BP is as accurate as RW, while again RWR performs better, achieving $100\%$ of success at $m=20$. 

The price of the improvement obtained by RW and RWR is the increased number of iterations. Future work will consider the development of faster strategies to get the desired local minimum. However, ADMM iterations are low-complex and the time spent for these experiments is acceptable ($10^5$ iterations require less than 20 seconds on a CPU @ 1.80GHz, RAM 16Gb).
\section{Conclusion}
In this paper, we have introduced a new  efficient framework to tackle the recovery of sparse binary signals acquired according to the compressed sensing paradigm. We have formulated the problem as local minimization of a non-convex, polynomial cost functional, which has a (local) minimum that corresponds to the desired conditions under very mild conditions. The search of such local minimum can be efficiently performed via iterative algorithms, such as reweighting procedures. In future work, we will study the conditions under which the desired minimum is the global minimum of the functional, and we will extend the approach to larger, non-binary alphabets and to systems with noise.


\begin{thebibliography}{10}
\providecommand{\url}[1]{#1}
\csname url@samestyle\endcsname
\providecommand{\newblock}{\relax}
\providecommand{\bibinfo}[2]{#2}
\providecommand{\BIBentrySTDinterwordspacing}{\spaceskip=0pt\relax}
\providecommand{\BIBentryALTinterwordstretchfactor}{4}
\providecommand{\BIBentryALTinterwordspacing}{\spaceskip=\fontdimen2\font plus
\BIBentryALTinterwordstretchfactor\fontdimen3\font minus
  \fontdimen4\font\relax}
\providecommand{\BIBforeignlanguage}[2]{{%
\expandafter\ifx\csname l@#1\endcsname\relax
\typeout{** WARNING: IEEEtran.bst: No hyphenation pattern has been}%
\typeout{** loaded for the language `#1'. Using the pattern for}%
\typeout{** the default language instead.}%
\else
\language=\csname l@#1\endcsname
\fi
#2}}
\providecommand{\BIBdecl}{\relax}
\BIBdecl

\bibitem{nak12}
U.~Nakarmi and N.~Rahnavard, ``{BCS}: Compressive sensing for binary sparse
  signals,'' in \emph{IEEE Military Communications Conference (MILCOM)}, 2012,
  pp. 1--5.

\bibitem{don06}
D.~L. Donoho, ``Compressed sensing,'' \emph{IEEE Trans. Inf. Theory}, vol.~52,
  no.~4, pp. 1289--1306, 2006.

\bibitem{kei17}
S.~Keiper, G.~Kutyniok, D.~G. Lee, and G.~E. Pfander, ``Compressed sensing for
  finite-valued signals,'' \emph{Linear Algebra and its Applications}, vol.
  532, no. Supplement C, pp. 570--613, 2017.

\bibitem{fou13}
S.~Foucart and H.~Rauhut, \emph{A Mathematical Introduction to Compressive
  Sensing}.\hskip 1em plus 0.5em minus 0.4em\relax New York: Springer, 2013.

\bibitem{sto10}
M.~Stojnic, ``Recovery thresholds for $\ell_1$ optimization in binary
  compressed sensing,'' in \emph{IEEE International Symposium on Information
  Theory (ISIT)}, 2010, pp. 1593--1597.

\bibitem{tia09}
Z.~Tian, G.~Leus, and V.~Lottici, ``Detection of sparse signals under
  finite-alphabet constraints,'' in \emph{IEEE International Conference on
  Acoustics, Speech and Signal Processing (ICASSP)}, 2009, pp. 2349--2352.

\bibitem{shi15}
M.~Shirvanimoghaddam, Y.~Li, B.~Vucetic, J.~Yuan, and P.~Zhang, ``Binary
  compressive sensing via analog fountain coding,'' \emph{IEEE Trans. Signal
  Process.}, vol.~63, no.~24, pp. 6540--6552, 2015.

\bibitem{ahn16}
J.~H. Ahn, ``Compressive sensing and recovery for binary images,'' \emph{IEEE
  Trans. Image Process.}, vol.~25, no.~10, pp. 4796--4802, 2016.

\bibitem{lee16}
N.~Lee, ``Map support detection for greedy sparse signal recovery algorithms in
  compressive sensing,'' \emph{IEEE Trans. Signal Process.}, vol.~64, no.~19,
  pp. 4987--4999, 2016, support detection method, and a sufficient condition
  for perfect signal recovery is derived for the case when the sparse signal is
  binary.

\bibitem{liu18}
\BIBentryALTinterwordspacing
T.~Liu and D.~G. Lee, ``Fast binary compressive sensing via $\ell_0$ gradient
  descent,'' 2018, https://arxiv.org/pdf/1801.09937.pdf. [Online]. Available:
  \url{https://arxiv.org/pdf/1801.09937.pdf}
\BIBentrySTDinterwordspacing

\bibitem{bay15}
A.~Bay, D.~Carrera, S.~M. Fosson, P.~Fragneto, M.~Grella, C.~Ravazzi, and
  E.~Magli, ``Block-sparsity-based localization in wireless sensor networks,''
  \emph{EURASIP Journal on Wireless Communications and Networking}, vol. 2015,
  no. 182, pp. 1--15, 2015.

\bibitem{fox13}
S.~M. Fosson, ``Binary input reconstruction for linear systems: A performance
  analysis,'' \emph{Nonlinear Analysis: Hybrid Systems}, vol.~7, no.~1, pp. 54
  -- 67, 2013.

\bibitem{fox15}
F.~Fagnani and S.~M. Fosson, ``Analysis of reduced-search bcjr algorithms for
  input estimation in a jump linear system,'' \emph{Signal Processing}, vol.
  108, pp. 341 -- 350, 2015.

\bibitem{tib96}
R.~Tibshirani, ``Regression shrinkage and selection via the lasso,''
  \emph{Journal of the Royal Statistical Society, Series B}, vol.~58, pp.
  267--288, 1996.

\bibitem{zha10MCP}
C.-H. Zhang, ``Nearly unbiased variable selection,'' \emph{Ann. Statist.},
  vol.~38, no.~2, pp. 894--942, 2010.

\bibitem{zha12}
C.-H. Zhang and T.~Zhang, ``A general theory of concave regularization for
  high-dimensional sparse estimation problems,'' \emph{Statist. Sci.}, vol.~27,
  no.~4, pp. 576 --593, 2012.

\bibitem{woo16}
J.~Woodworth and R.~Chartrand, ``Compressed sensing recovery via nonconvex
  shrinkage penalties,'' \emph{Inverse Problems}, vol.~32, no.~7, pp.
  75\,004--75\,028, 2016.

\bibitem{fox16}
S.~M. Fosson, J.~Matamoros, C.~Ant\'{o}n-Haro, and E.~Magli, ``Distributed
  recovery of jointly sparse signals under communication constraints,''
  \emph{IEEE Trans. Signal Process.}, vol.~64, no.~13, pp. 3470--3482, 2016.

\bibitem{can08rew}
E.~J. Cand\`es, M.~B. Wakin, and S.~Boyd, ``Enhancing sparsity by reweighted
  $\ell_1$ minimization,'' \emph{Journ. Fourier Anal. Appl.}, vol.~14, no. 5-6,
  pp. 877--905, 2008.

\bibitem{lasbook}
J.-B. Lasserre, \emph{An introduction to polynomial and semi-algebraic
  optimization}, ser. Cambridge texts in applied mathematics.\hskip 1em plus
  0.5em minus 0.4em\relax Cambridge, UK: Cambridge University Press,, 2015.

\bibitem{tib13}
R.~J. Tibshirani, ``{The {L}asso problem and uniqueness},'' \emph{Electronic
  Journal of Statistics}, vol.~7, pp. 1456--1490, 2013.

\bibitem{zou08_LLA}
H.~Zou and R.~Li, ``One-step sparse estimates in nonconcave penalized
  likelihood models,'' \emph{Annals of Statistics}, vol.~36, no.~4, p. 1509,
  2008.

\bibitem{faz03}
M.~Fazel, H.~Hindi, and S.~Boyd, ``Log-det heuristic for matrix rank
  minimization with applications to {H}ankel and {E}uclidean distance
  matrices,'' in \emph{IEEE Proc. American Control Conference (ACC)}, vol.~3,
  2003, pp. 2156--2162.

\bibitem{fox18}
S.~M. Fosson, ``A biconvex analysis for lasso $\ell_1$ reweighting,''
  \emph{IEEE Signal Process. Lett.}, vol. early access, no.~nn, pp. 1--1, 2018.

\bibitem{asi13}
M.~Asif and J.~Romberg, ``Fast and accurate algorithms for re-weighted $\ell
  _{1} $-norm minimization,'' \emph{IEEE Trans. Signal Process.}, vol.~61,
  no.~23, pp. 5905--5916, 2013.

\bibitem{ahm15}
R.~Ahmad and P.~Schniter, ``Iteratively reweighted $\ell_1$ approaches to
  sparse composite regularization,'' \emph{IEEE Trans. Computat. Imag.},
  vol.~1, no.~4, pp. 220--235, 2015.

\bibitem{cal16}
M.~Calvo-Fullana, J.~Matamoros, C.~Ant\'{o}n-Haro, and S.~M. Fosson,
  ``Sparsity-promoting sensor selection with energy harvesting constraints,''
  in \emph{IEEE International Conference on Acoustics, Speech and Signal
  Processing (ICASSP)}, 2016, pp. 3766--3770.

\bibitem{fox14}
S.~M. Fosson, J.~Matamoros, C.~Ant\'{o}n-Haro, and E.~Magli, ``Distributed
  support detection of jointly sparse signals,'' in \emph{IEEE International
  Conference on Acoustics, Speech and Signal Processing (ICASSP)}.\hskip 1em
  plus 0.5em minus 0.4em\relax IEEE, 2014, pp. 6434--6438.

\bibitem{boy10}
S.~Boyd, N.~Parikh, E.~Chu, B.~Peleato, and J.~Eckstein, ``Distributed
  optimization and statistical learning via the alternating direction method of
  multipliers,'' \emph{Found. Trends Mach. Learn.}, vol.~3, no.~1, pp. 1 --
  122, 2010.

\end{thebibliography}
\end{document}